\documentclass[11pt]{article}
\usepackage{latexsym,amssymb}
\usepackage{amsfonts}

\def\demo{\noindent{\bf Proof. }}

\def\QED{\hfill$\Box$}

\newtheorem{Theorem}{Theorem}[section]
\newtheorem{Lemma}[Theorem]{Lemma}
\newtheorem{Corollary}[Theorem]{Corollary}
\newtheorem{Proposition}[Theorem]{Proposition}
\newtheorem{Remark}[Theorem]{Remark}
\newtheorem{Example}[Theorem]{Example}

\newtheorem{Definition}[Theorem]{Definition}

\begin{document}

\title{Lectures on quasi-invariants of Coxeter groups and the Cherednik algebra}
\author{Pavel Etingof\\
Massachusetts Institute of Technology\\ Department of
Mathematics\\ Cambridge, MA 02139, USA\\ {\rm email:
etingof@math.mit.edu} \and Elisabetta Strickland\\ 
Dipartimento di Matematica\\ Universit\'a di Roma
"Tor Vergata"\\ Via della Ricerca Scientifica - 00133 Roma, Italy.
\\ {\rm email: strickla@mat.uniroma2.it} }
\maketitle

\centerline{\bf Introduction}

This paper arose from a series of three lectures given by the first author
at Universit\'a di Roma ``Tor Vergata'' in January 2002, when the
second author extended and improved her notes of these lectures. 
It contains an elementary introduction for non-specialists to the
theory of quasi-invariants (but no original results).

Our main object of study is the variety
$X_m$ of quasi-invariants for a finite Coxeter group. 
This very interesting singular algebraic variety 
arose in a work of O.Chalykh and A.Veselov about 10 years ago, 
as the spectral variety of the quantum Calogero-Moser system. 
We will see that despite being singular, this variety has very nice
properties (Cohen-Macaulay, 
Gorenstein, simplicity of the ring of differential operators,
explicitly given Hilbert series). It is interesting that although
the definition of $X_m$ is completely elementary, 
to understand the geometry of $X_m$ it is helpful to use representation 
theory of the rational degeneration of Cherednik's double affine
Hecke algebra, and the theory of integrable systems. 
Thus, the study of $X_m$ leads us to a junction of
three subjects -- integrable systems, representation theory, and 
algebraic geometry. 

The content of the paper is as follows. 
In Lecture 1 we define the ring of quasi-invariants 
for a Coxeter group, and discuss its elementary 
properties (with proofs), as well as deeper properties, such as
Cohen-Macaulay, Gorenstein property, and the Hilbert series
(whose proofs are partially postponed until Lecture 3).
In Lecture 2, we explain the origin of the ring of
quasi-invariants in the theory of integrable systems, and 
introduce some tools from integrable systems, such as the
Baker-Akhieser function. Finally, in Lecture 3, we 
 develop the theory of the rational Cherednik algebra,
the representation-theoretic techniques due to Opdam and Rouquier,
and finish the proofs of the geometric statements from Chapter 1.

{\bf Acknowledgments.} The authors are grateful to Corrado De
Concini for useful discussions. The work of P.E. was partially
supported by his NSF grant DMS-9988796 and was done in part for
the Clay Mathematics Institute. The series of lectures in ``Tor Vergata" has been finaced by the MIUR
"Progetto Azioni di Gruppi" " of which E.S. is a member.             

\section{Lecture 1}
\subsection {Definition of quasi-invariants}

In this lecture we will define the ring of quasi-invariants $Q_m$
and discuss its main properties. 

We will work over the field of 
complex numbers $\Bbb C$. Let $W$ be a finite Coxeter group i.e. a 
finite group generated by reflections.  Let  us denote by $\frak h$ 
its reflection representation.  A typical example
is the Weyl group of a simple Lie algebra acting on a 
Cartan subalgebra $\frak h$. In the case the Lie algebra is 
$s\ell(n)$, we have that $W$ is the symmetric group $S_n$ on $n$ 
letters and $\frak h$ is the space of diagonal traceless $n\times n$ 
matrices. 


Let $\Sigma\subset W$ denote the set 
of reflections. Clearly, $W$ acts on $\Sigma$ by conjugation. Let 
$m:\Sigma\to \Bbb Z_+$ be a function on $\Sigma$ taking non negative 
integer values, which is $W$-invariant. 
The number of orbits of $W$ on $\Sigma$ is generally very small.
For example, if $W$ 
is the Weyl group of a simple Lie algebra of ADE type, then 
$W$ acts transitively on $\Sigma$, so $m$ is a constant
function. 

For each reflection $s\in \Sigma$, choose 
$\alpha_s\in \frak h^*-\{0\}$ so that 
$\alpha_s(sx)=-\alpha_s(x)$ (this means that  the hyperplane
given by the 
equation $\alpha_s=0$ is the reflection hyperplane for
$s$).

\begin 
{Definition} \cite{CV1,CV2}
{\rm A polynomial $q\in \Bbb C[\frak h]$ is said {\it 
$m$-quasiinvariant} with respect to $W$, if for any $s\in
\Sigma$, the polynomial
$q(x)-q(sx)$ is divisible by $\alpha_s(x)^{2m_s+1}$.}
\end {Definition} 

We will denote by $Q_m$ the space of 
$m$-quasiinvariant  polynomials for $W$.

Notice that every element 
of  $\Bbb C[\frak h]$ is a $0$-quasiinvariant, and that every $W$ 
invariant is an $m$-quasiinvariant for any $m$. Indeed if $q\in \Bbb 
C[\frak h]^W$, then we have $q(x)-q(sx)=0$ for all $s\in \Sigma$, and $0$ is 
divisible by all powers of $\alpha_s(x)$. Thus in a way, $\Bbb 
C[\frak h]^W$ can be viewed as the set of $\infty$-quasinvariants.

\begin{Example}{\rm   The group $W=\Bbb Z/2$ acts on $\frak h=\Bbb C$ by 
$s(v)=-v$.  In this case $m\in \Bbb Z_+$ and $\Sigma=\{s\}$. So this 
definition says that $q$ is in $Q_m$ iff $q(x)-q(-x)$ is divisible 
by $x^{2m+1}$. It is very easy to write a basis  of $Q_m$. It is 
given by the polynomials
$\{x^{2i}|i\geq 0\}\cup \{x^{2i+1}|i\geq 
m\}$.
}\end{Example}

\subsection{Elementary properties of $Q_m$.}

Some elementary properties of $Q_m$ are collected in the
following proposition. 

\begin 
{Proposition}\label{prima} (See \cite{FV} and references therein).

1) $\Bbb C[\frak h]^W\subset Q_m \subseteq 
\Bbb C[\frak h],$
$Q_0=\Bbb C[\frak h]$, $Q_m\subset Q_{m'}$ if 
$m\ge m'$, $\cap_mQ_m=\Bbb C[\frak h]^W$.

2) $Q_m$ is a subring of $\Bbb 
C[\frak h]$.

3) The fraction field of $Q_m$  is equal to  $\Bbb 
C(\frak h)$.

4) $Q_m$ is a finite $\Bbb C[\frak h]^W$-module and a 
finitely generated algebra.
$\Bbb C[\frak h]$ is a finite 
$Q_m$-module.
\end{Proposition}
\demo 1) is immediate and  has partly 
been remarked already.

2) Clearly $Q_m$ is closed under sum. Let 
$p,q\in Q_m$. Let $s\in \Sigma$. Then 
$$p(x)q(x)-p(sx)q(sx)=p(x)q(x)-p(sx)q(x)+p(sx)q(x)-p(sx)q(sx)=$$ 
$$=(p(x)-p(sx))q(x)+p(sx)(q(x)-q(sx))$$
Since both $p(x)-p(sx)$ and 
$q(x)-q(sx)$ are divisible by $\alpha_s^{2m_s+1}$, we deduce that  
$p(x)q(x)-p(sx)q(sx)$ is also divisible by $\alpha_s^{2m_s+1}$, 
proving the claim.

3) Consider the polynomial 

$$\delta_{2m+1}(x)=\prod_{s\in \Sigma}\alpha_s(x)^{2m_s+1}$$
This 
polynomial is uniquely defined up to scaling.  One has 
$\delta_{2m+1}(sx)=-\delta_{2m+1}(x)$. Take $f(x)\in \Bbb C[\frak 
h]$. We claim that
$f(x)\delta_{2m+1}(x)\in Q_m$. As a matter of 
fact, 
$$f(x)\delta_{2m+1}(x)-f(sx)\delta_{2m+1}(sx)=(f(x)+f(sx))\delta_{2m+1}(x).$$
and by its definition $\delta_{2m+1}(x)$ is divisible by 
$\alpha_s(x)^{2m_s+1}$ for all $s\in \Sigma$. This implies 3).

4) By 
Hillbert's theorem on the finiteness of invariants, we get that $\Bbb 
C[\frak h]^W$ is a finitely generated algebra over $\Bbb C$ and $\Bbb 
C[\frak h]$ is a finite $\Bbb C[\frak h]^W$-module and hence a 
finite $Q_m$-module, proving the second part of 4).

Now $Q_m\subset 
\Bbb C[\frak h]$ is  a submodule of the finite module $\Bbb C[\frak 
h]$ over the Noetherian ring $\Bbb C[\frak h]^W$. Hence it is finite. 
This immediately implies that $Q_m$ is a finitely generated 
algebra. $\square$

{\bf Remark.} 
In fact, since $W$ is a finite Coxeter
group, a celebrated  result 
of Chevalley says that the algebra $\Bbb C[\frak h]^W$ is not only a finitely 
generated $\Bbb C$-algebra but actually a free (=polynomial)
algebra.
Namely, it has the form $\Bbb C[q_1, \ldots 
, q_n]$, where the $q_i$ are homogeneous polynomials of some degrees 
$d_i$. Furthermore, if we 
denote by $H$ the subspace of $\Bbb C[\frak h]$ of harmonic 
polynomials, i.e. of polynomials killed by $W$ invariant differential 
operators with constant coefficients without constant term, then the 
multiplication map 
$$\Bbb C[\frak h]^W\otimes H\to \Bbb C[\frak h]$$
is an isomorphism of $\Bbb C[\frak h]^W$-  and of $W$-modules. In 
particular, $\Bbb C[\frak h]$ is a free $\Bbb C[\frak h]^W$ module of 
rank $|W|$.

\subsection{The variety $X_m$ and its bijective normalization}

Using Proposition~\ref{prima}, we can define the 
irreducible affine variety $X_m=Spec(Q_m)$.
The inclusion $Q_m\subset 
\Bbb C[\frak h]$ induces a morphism
$$\pi:\frak h\to X_m$$
which 
again by   Proposition~\ref{prima}, is birational and surjective. 
(Notice that in particular this implies that $X_m$ is singular for all 
$m\ne 0$).

In fact, not only is $\pi$  birational, but  a stronger 
result is true.
\begin {Proposition} (Berest, see \cite{BEG})\label{seconda}
$\pi$ is a 
bijection.
\end{Proposition}\demo  By the above remarks, we only have 
to show that $\pi$ is injective. In order to achieve this, we need to 
prove that quasi-invariants separate points of $\frak h$, i.e. if 
$z,y\in \frak h$ and $z\neq y$, then there exist $p\in Q_m$ such 
that $p(z)\neq p(y)$. This is obtained in the following way.  Let 
$W_z\subset W$ be the stabilizer of $z$ and choose $f\in \Bbb C[\frak 
h]$ such that $f(z)\neq 0$, $f(y)=0$. Set 
$$p(x)=\prod_{s\in \Sigma 
, sz\neq z}\alpha_s(x)^{2m_s+1}\prod_{w\in W_z}f(wx).$$
We claim that 
$p(x)\in Q_m$. Indeed, let $s\in \Sigma$ and assume that $s(z)\neq z$.

We have by definition $p(x)=\alpha_s(x)^{2m_s+1}\tilde p(x)$, with 
$\tilde p(x)$ a polynomial. So

$$p(x)-p(sx)=\alpha_s(x)^{2m_s+1}\tilde 
p(x)-\alpha_s(sx)^{2m_s+1}\tilde p(sx)=\alpha_s(x)^{2m_s+1}(\tilde 
p(x)+\tilde p(sx))$$
If on the other hand, $sz=z$, i.e. $s\in W_z$, 
then $s$ preserves the set $W\setminus W_z$, 
and hence preserves 
$\prod_{s\in W\setminus W_z}\alpha_s(x)^{2m_s+1}$
(as it acts by $-1$ on the products of the same terms both over
$W$ and over $W_z$). Since 
$\prod_{w\in W_z}f(wx)$ 
is $W_z$ invariant, we deduce that $p(x)-p(sx)=0$, so that in 
this case
$p(x)-p(sx)$ also is divisible by $\alpha_s(x)^{2m_s+1}$.

To 
finish, notice that $p(z)\neq 0$. Indeed, for a reflection $s$, 
$\alpha_s$ vanishes exactly on the fixed points of $s$, so that 
$\prod_{s\in \Sigma , sz\neq z}\alpha_s(z)^{2m_s+1}\neq 0$. Also for 
all  $w\in W_z$ $f(wz)=f(z)\neq 0$.  On the other hand, it is clear 
that $p(y)=0$.\QED

\begin{Example}{\rm  $W=\Bbb Z/2$. As we have 
already seen, $Q_m$ has a basis given by the monomials 
$\{x^{2i}|i\geq 0\}\cup \{x^{2i+1}|i\geq m\}$. From this we deduce 
that setting $z=x^2$ and $y=x^{2m+1}$, $Q_m=\Bbb 
C[y,z]/(y^2-z^{2m+1})=\Bbb C[K]$, where $K$ is the plane curve with a 
cusp at the origin, given by the equation $y^2=z^{2m+1}$. The map
$\pi :\Bbb 
C\to K$ is given by $\pi(t)=(t^2,t^{2m+1})$ which is clearly 
bijective.
}\end{Example}

\subsection{Further properties of $X_m$}

Let us get to some deeper properties of 
quasi-invariants. Let $X$ be an irreducible affine variety over $\Bbb 
C$ and $A=\Bbb C[X]$. Recall that,  by Noether normalization Lemma, 
there exists $f_1,\ldots ,f_n\in \Bbb C[X]$ which are algebraically 
independent over $\Bbb C$ and such that $\Bbb C[X]$ is a finite 
module over the polynomial ring $\Bbb C[f_1,\ldots ,f_n]$. This means 
that we have a finite morphism of $X$ onto an affine space. 

\begin 
{Definition} {\rm $A$ (and $X$) is said to be 
Cohen-Macaulay if there exist $f_1,\ldots 
f_n$ as above, with the property that $\Bbb C[X]$ is a locally free 
module over $\Bbb C[f_1,\ldots ,f_n]$. (Notice that by the 
Quillen-Suslin theorem, this is equivalent to saying that $A$ is
a free module).}
\end {Definition}

{\bf Remark.} If $A$ is Cohen-Macaulay, then for any 
$f_1,\ldots ,f_n$ which are algebraically 
independent over $\Bbb C$ and such that $\Bbb C[X]$ is a finite 
module over the polynomial ring $\Bbb C[f_1,\ldots ,f_n]$, $A$ 
is a locally free $\Bbb C[f_1,\ldots ,f_n]$-module,
\cite{Eis}.

\begin {Theorem}\label{pip}(\cite{EG2},\cite{BEG}, 
conjectured in \cite{FV})  $Q_m$ is Cohen-Macaulay.
\end {Theorem}
 
Notice that, using Chevalley's result that $\Bbb C[\frak h]^W$ is a 
polynomial ring, in order to prove Theorem~\ref{pip}  it will suffice 
to prove:
\begin {Theorem} ({\rm \cite{EG2,BEG}, conjectured in \cite{FV}})
  \label{freeness} $Q_m$ is a free
$\Bbb C[\frak h]^W$ module.
\end
{Theorem} 

A proof of this theorem will be given at the end of Lecture 3.
This proof follows \cite{BEG} (the original proof of \cite{EG2}
is shorter but somewhat less conceptual). 

\subsection{The Poincar\'e series of $Q_m$}

Consider now the Poincar\'e series

$$h_{Q_m}(t)=\sum_{r\geq 0} {\rm dim}Q_m[r]t^r$$
$Q_m[r]$ denoting 
the graded component of $Q_m$ of degree $r$.

For every irreducible 
representation $\tau\in \widehat W$, define

$$\chi_{\tau}(t)=\sum_{r\geq 0}{\rm dim Hom}_W(\tau,\Bbb C[\frak 
h][r])t^r.$$
Consider the element in the group ring $\Bbb Z[W]$

$$\mu_m=\sum_{s\in \Sigma}m_s(1-s).$$
The $W$ invariance of $m$ 
implies that $\mu_m$ lies in the center of $\Bbb Z[W]$.
Hence it is 
clear that $\mu$ acts as a scalar, $\xi_m(\tau)$,  on $\tau$. 
\begin 
{Lemma} The scalar $\xi_m(\tau)$ is an integer.
\end {Lemma} 
\demo  
$\Bbb Z[W]$ and hence also its center, is a finite $\Bbb Z$-module. 
This clearly implies that  
$\xi_m(\tau)$ is an algebraic integer. 
Thus to prove that $\xi_m(\tau)$ is an integer, it suffices to see 
that $\xi_m(\tau)$ is a rational number. Set $d_{\tau}$ equal to the 
degree of $\tau$ and $d_{\tau,s}$ equal to the dimension of space of 
$s$ invariants in $\tau$. Taking traces we get

$$d_{\tau}\xi_m(\tau)=\sum_{s\in \Sigma}2m_s(d_{\tau}-d_{\tau,s})$$

which gives the rationality of  $\xi_m(\tau)$.\QED

One has:
\begin 
{Theorem} 
\begin{equation}\label{unouno} h_{Q_m}(t)=\sum_{\tau\in\widehat 
W}d_{\tau}t^{\xi_m(\tau)}\chi_{\tau}(t)\end{equation}
\end{Theorem}

{\bf Remark.} This theorem was proved in \cite{FeV} 
modulo Theorem \ref{pip} (conjectured in \cite{FV}) 
using the Matsuo-Cherednik correspondence. Thus, 
Theorem \ref{unouno} follows from \cite{FeV} and \cite{EG2}. 
Another proof of this theorem is given in \cite{BEG}; this is
the proof we will discuss below. 

\begin{Example}{\rm  If $m=0$, since $Q_0=\Bbb C[\frak h]$, 
the theorem says that
$$h_{Q_0}(t)={1\over (1-t)^n}=\sum_{\tau\in\widehat 
W}d_{\tau}\chi_{\tau}(t)$$
Indeed, as a $W$-module one has
$$\Bbb 
C[\frak h]=\oplus_\tau\tau\otimes {\rm Hom}_W(\tau,\Bbb C[\frak h]).$$

}\end{Example}
\begin{Example}{\rm  If $W=\Bbb Z/2$, then $\widehat 
W=\{+, -\}$, $+$ (respectively $-$) denoting the trivial 
(respectively the sign) representation. One has
$$\Bbb C[x]=\Bbb C[x^2]\oplus 
\Bbb C[x^2]x$$
where $\Bbb C[x^2]=\Bbb C[x]^W$ and $\Bbb C[x^2]x$ is 
the isotypic component of the sign representation. Thus
$$ 
\chi_{+}(t)={1\over 1-t^2},\ \  \chi_-(t)={t\over 1-t^2},$$
$\mu_m=m(1-s).$ Thus $\xi_m(+)=0$, $\xi_m(-)=2m$.
We deduce that

$$h_{Q_m}(t)={1\over 1-t^2}+{t^{2m+1}\over 1-t^2}$$
as we already 
know.
}\end{Example}
Recall now that as a graded $W$-module $\Bbb 
C[\frak h]$ is isomorphic to 
$\Bbb C[\frak h]^W\otimes H$, $H$ being 
the space of harmonic polynomials. We deduce that the $\tau$-isotypic 
component in $\Bbb C[\frak h]$ is isomorphic to $\Bbb C[\frak 
h]^W\otimes H_{\tau}$. Set $K_{\tau}(t)=\sum_{r\geq 0}{\rm dimHom}(\tau 
,H[r])t^r$. This is a polynomial called the Kostka polynomial 
relative to $\tau$. We deduce that
\begin{equation}\label{unodue} 
\chi_{\tau}(t)={K_{\tau}(t)\over 
\prod_{i=1}^{n}(1-t^{d_i})}\end{equation}
Also, if $\tau'=\tau\otimes 
\varepsilon$, $\varepsilon$ being the sign representation, one has

$$K_{\tau'}(t)=K_{\tau}(t^{-1})t^{|\Sigma|}$$
Set now 

$$P_m(t)=\sum_{\tau\in \widehat 
W}d_{\tau}t^{\xi_m(\tau)}K_{\tau}(t)$$
We have
\begin 
{Proposition}\cite{FeV} \label{terza} $$h_{Q_m}(t)={P_m(t)\over 
\prod_{i=1}^{n}(1-t^{d_i})}.$$
Furthermore 
$P_m(t)=t^{\xi_m(\varepsilon)+|\Sigma|}P_m(t^{-1})$.
\end 
{Proposition}
\demo Substituting the expression (\ref{unodue})  for 
$\chi_{\tau}(t)$ in (\ref{unouno}) and using the definition of 
$P_m(t)$, we get
$$h_{Q_m}(t)=\sum_{\tau\in\widehat 
W}d_{\tau}t^{\xi_m(\tau)}{K_{\tau}(t)\over 
\prod_{i=1}^{n}(1-t^{d_i})}={P_m(t)\over 
\prod_{i=1}^{n}(1-t^{d_i})}$$
as desired.

Now notice that 
$$\xi_m(\tau)+\xi_m(\tau')=\sum_{s\in 
\Sigma}2m_s=\xi_m(\varepsilon)$$
Using this we get

$$t^{\xi_m(\varepsilon)+|\Sigma|}P_m(t^{-1})=\sum_{\tau\in \widehat 
W}d_{\tau}t^{\xi_m(\varepsilon)-\xi_m(\tau)}t^{|\Sigma|}K_{\tau}(t^{-1})=\sum_{\tau'\in 
\widehat W}
d_{\tau'}t^{\xi_m(\tau')}K_{\tau'}(t)=P_m(t)$$
as 
desired.\QED

From this we deduce
\begin {Theorem}\label{goren}  ({\rm 
\cite{EG2,BEG,FeV}, conjectured in \cite{FV}}) The ring $Q_m$ of 
$m$-quasi-invariants is Gorenstein.
\end {Theorem}  
\demo By Stanley's theorem (see \cite{Eis}), a positively graded 
Cohen-Macaulay domain $A$ is Gorenstein iff its Poincare series 
is a rational function $h(t)$ satisfying the equation $h(t^{-1})=
(-1)^nt^lh(t)$, where $l$ is an integer and $n$ the 
dimension of the spectrum of $A$.
Thus the result 
follows immediately from Proposition~\ref{terza}.\QED

\subsection{The ring of differential operators on $X_m$}

Finally, let us introduce the ring  ${\cal D}(X_m)$  of 
differential operators on $X_m$, that is the ring of differential 
operators with coefficients in $\Bbb C(\frak h)$ mapping $Q_m$ to 
$Q_m$. 
It is clear that this definition coincides with the well known Grothendieck's definition. 

\begin {Theorem} \cite{BEG} \label{Dsimple} ${\cal 
D}(X_m)$ is a simple algebra.
\end{Theorem}
\begin {Remark} {\rm a) 
The ring of differential operators on a smooth affine algebraic 
variety is always simple.

b) By a result of M. van den Bergh \cite{VdB}, for a non-smooth variety, 
the simplicity of the ring of  differential operators 
implies the Cohen-Macaulay property of this variety.}
\end{Remark}

\section {Lecture 2}
We will now see how the ring $Q_m$ appears in the 
theory of completely integrable systems. 

\subsection{Hamiltonian mechanics and integrable systems}

Recall the basic setup of Hamiltonian mechanics \cite{Ar}. 
Consider a mechanical system with configuration space $X$ (a
smooth manifold). Then the phase space of this system is $T^*X$,
the cotangent bundle on $X$. The space $T^*X$ is naturally a
symplectic manifold, and in particular we have an operation of
Poisson bracket on functions on $T^*X$. A point of $T^*X$ is 
a pair $(x,p)$, where $x\in X$ is the position and $p\in T^*_xX$
is the momentum. Such pairs are called states of the system. 
The dynamics of the system $x=x(t)$, $p=p(t)$
depends on the Hamiltonian, or energy function, $E(x,p)$ on
$T^*X$. Given $E$ and the initial state $x(0),p(0)$, one can
recover the dynamics $x=x(t),p=p(t)$ from Hamilton's differential
equations $\frac{df(x,p)}{dt}=\lbrace{E,f\rbrace}$. If
$X$ is locally identified with $\Bbb R^n$ by choosing
coordinates $x_1,...,x_n$, then $T^*X$ is locally identified with 
$\Bbb R^{2n}$ with coordinates $x_1,...,x_n,p_1,...,p_n$. 
In these coordinates, Hamilton's equations may be written in
their standard form
$$
\dot 
x_i={\partial E\over \partial p_i},\  \ \  \dot p_i=-{\partial E\over 
\partial x_i}.
$$

A function $I(x,p)$ is called an integral of motion for our 
system  if $\{I,E\}=0$. Integrals of motion are useful, since 
for any such integral $I$ the function $I(x(t),p(t))$ 
is constant, which allows one to reduce the number of variables by
$2$. Thus, if we are given $n$ functionally independent integrals of motion
$I_1, \ldots I_n$ with $\{I_l,I_k\}=0$ for all $1\leq l,k\leq n$,
then all $2n$ variables $x_i,p_i$ can be excluded, and the system
can be completely solved by quadratures. Such situation is called complete (or
Liouville) integrability. 

\subsection{Classical Calogero-Moser system}

Quasi-invariants are related to many-particle systems. 
Consider a system 
of $n$ particles on the real line $\Bbb R$. 
A potential is an even function
$$U(x)=U(-x).$$
Two particles at points 
$a,b$ have energy of interaction $U(a-b)$. The total energy of our system of particles 
is
$$E=\sum_{i=1}^n{p_i^2\over 2}+\sum_{i<j}U(x_i-x_j).$$
Here, $x_i$ are the coordinates of the particles, $p_i$ their
momenta.  
The dynamics of the particles $x_i=x_i(t)$, $p_i=p_i(t)$
is governed by the Hamilton equations with energy function $E$. 

This is a system of nonlinear differential equations, which in general
may be difficult to solve explicitly. 
However, for special potentials this system 
may be completely integrable. For instance, we will see that 
it is so for the Calogero-Moser potential,
$$U(x)={\gamma\over x^2},$$
$\gamma$ being a 
constant.

The Calogero-Moser system has a generalization to 
arbitrary Coxeter groups. Namely, 
consider a finite  group $W$  generated by reflections 
acting on the space $\frak h$, and keep the notation of the
previous section. Fix a $W$-invariant nondegenerate 
scalar product $(-,-)$ on $\frak h$. It determines a scalar
product on $\frak h^*$. Define the ``energy function''
$$E(x,p)={(p,p)\over 2}+{1\over 
2}\sum_{s\in \Sigma}{\gamma_s(\alpha_s,\alpha_s)\over 
\alpha_s(x)^2}.$$
on $T^*\frak h=\frak h\times \frak h^*$, where 
$\gamma : \Sigma\to \Bbb C$ is   a $W$-invariant function. 
Notice that although $\alpha_s$
is defined up to a non zero constant, 
by homogeneity, $E$ is independent of the choice of $\alpha_s$. We 
will call the system defined by $E$ the Calogero-Moser system for $W$.

If 
$W$ is the symmetric group $S_n$, $\frak h= \Bbb C^n$, then
$\Sigma$ is the set of transpositions $s_{i,j}$, $i<j$ and we can 
take $\alpha_s=e_i-e_j$,  Then we clearly obtain the usual Calogero-Moser 
system. 

Below we will see that the Calogero-Moser system for $W$ is completely
integrable. 

\subsection{Quantum Calogero-Moser system}

Let us now discuss quantization of the Calogero-Moser system.
We start by quantizing the 
energy $E$ by formally making the substitution
$$p_j\Rightarrow 
i\hbar {\partial\over \partial x_j}$$
where $\hbar$ is a 
parameter (Planck constant). 
This yields the 
Schr\"odinger operator
$$\widehat E:=-{\hbar^2\over 2}\Delta+{1\over 
2}\sum_{s\in \Sigma}{\gamma_s(\alpha_s,\alpha_s)\over \alpha_s^2},$$
where $\Delta$ denotes the Laplacian.

In particular, in the case of 
$W=S_n$ we have
$$\widehat E=-{\hbar^2\over 2}\Delta+\sum_{i<j}{c\over 
(x_i-x_j)^2},$$
where $\Delta=\sum_i {\partial^2\over\partial x_i^2}$.

Setting $\beta_s=\frac{\gamma_s}{2\hbar^2}$, we will from 
now on consider the operator
$$H:=-{2\over \hbar^2}\widehat 
E=\Delta-\sum_{s\in \Sigma} {\beta_s(\alpha_s,\alpha_s)\over 
\alpha_s^2(x)}.$$

This operator is called the Calogero-Moser operator.

We want to study the stationary Schr\"odinger 
equation:
\begin{equation}
\label{due0}H\psi=\lambda \psi, \ \ \ \ \ \ \lambda\in \Bbb C.\end{equation}

Similarly to the classical case, 
for a general Schr\"odinger operator $H$, it is hard
to say anything explicit about solutions of this equation, 
but for the Calogero-Moser operator the situation is much
better. 

\begin {Definition} {\rm A quantum integral of $H$ is a differential 
operator $M$ such that $$[H,M]=0.$$}
\end {Definition} 

We are going 
to show that there are plenty of quantum integrals  of $H$, namely that
there are n commuting algebraically independent quantum integrals 
$M_1,\ldots ,M_n$  of $H$.  By definition, this means that the 
quantum Calogero-Moser system is completely integrable.

Once we have 
found $M_1,\ldots ,M_n$,   remark that for fixed constants 
$\mu_1,\ldots ,\mu_n$, the space of solutions of the system
$$\left 
\{\matrix {M_1\psi=\mu_1\psi \cr \cdots\cdots \cr M_n\psi=\mu_n 
\psi}\right .$$
is clearly stable under $H$. 
We will see that this space is in fact finite dimensional.
Therefore, the operators $M_i$ allow one to reduce 
solving the partial differential equation $H\psi=\lambda\psi$ 
to solving a system of ordinary linear differential equations. 
This phenomenon is called quantum complete integrability.

\subsection{The algebra of differential-reflection operators}
We are now going to 
explain how to find quantum integrals for $H$, using the Dunkl-Cherednik method.

First let us fix some notation. Given a smooth 
affine variety $X$, we will denote by ${\cal D}(X )$ the ring of 
differential operators on $X$ .   We are going to consider the case 
in which $X$ is the open set $U$ in $\frak h$ which is the complement 
of the divisor of the equation $\prod_{s\in \Sigma}\alpha_s(x)$.  Clearly 
${\cal D}(U)={\cal D}(\frak h)[1/\delta(x)]$. 

\begin{Lemma}\label{inva}
An element of ${\cal D}(U)$ is completely determined by its
action on $\Bbb C[U]^W=\Bbb C[U/W]$. 
\end{Lemma}

\demo
Recall that the 
quotient map $\pi: U\to U/W$ is finite and unramified. This implies 
that  $${\cal D}(U)={\Bbb C}[U]\otimes_{{\Bbb C}[U/W]}{\cal D}(U/W).$$
From this we obtain that if  $P\in {\cal D}(U)$ is such that $Pf=0$ for all 
$f\in \Bbb C[U/W]$ then $P=0$.
\QED

We also have the 
operators on $\Bbb C[U]$ given by the action of $W$. We will denote 
by $\cal A$ the algebra of operators on $U$ generated by ${\cal D}(U)$ 
and $W$.
We have:

\begin {Proposition} ${\cal A}={\cal D}(U)\rtimes W$ 
i.e. every element in $A\in \cal A$ can be uniquely written as a 
linear combination
$$A=\sum_{w\in W}P_ww$$
with $P_w\in {\cal D}(U)$.
\end{Proposition}

\demo 
The fact that every element in $\cal A$ can 
be expressed as a linear combination $\sum_{w\in W}P_ww$ is
clear. To show that such an expression is unique,
assume $\sum_{w\in W}P_ww=0$. Take 
$f\in \Bbb C[U]$ such that ${^w}f\neq {^{u}}f$ for all $w\neq u$ in 
$W$. Then

$$\sum_{w\in W}P_w{^w\negthinspace f}^iw=\sum_{w\in W}P_wwf^i=0$$
for all 
$i\geq 0$. 
We deduce that given $g\in \Bbb C[U/W]$ one has
$$(\sum_{w\in 
W}P_w{^w\negthinspace f}^i)(g)=0.$$
Thus by Lemma \ref{inva}, $\sum_{w\in W}P_w{^w\negthinspace f}^i=0$ for all $i$. 
Therefore,
$P_w\prod_{w\neq u}({^w}f- 
{^{u}}f)=0$ and hence $P_w=0$, for all $w\in W$, as desired.
\QED

Take $A\in \cal  A$ and write
$$A=\sum_{w\in W}P_ww.$$
We set 
$m(A)=\sum_{w\in W}P_w\in {\cal D}(U).$
Notice that if $f$ is a $W$-invariant function, 
then clearly $A(f)=m(A)(f)$ and that, by what we 
have seen in the proof of the above proposition, $m(A)$ is completely 
determined by its action on invariant functions. 

In general, $m$ is not a 
homomorphism. However:
\begin {Proposition} Let ${\cal A}^W\subset \cal 
A$ denote the subalgebra of $W$-invariant elements.
Then the 
restriction of $m$ to ${\cal A}^W$ is an algebra 
homomorphism.
\end{Proposition}
\demo If $A\in {\cal  A}^W$, then 
clearly $m(A)$ is $W$-invariant. Now if we take $A,B\in {\cal 
A}^W$ and $f$ a $W$-invariant function we have that $B(f)$ is
also $W$-invariant. So
$$m(AB)(f)=(AB)(f)=A(B(f))=A(m(B)(f))=m(A)(m(B)(f)).$$

Thus $m(AB)$ and $m(A)m(B)$ coincide on $W$-invariant functions and 
hence coincide.\QED

\subsection{Dunkl operators and symmetric quantum integrals}

Fix a $W$ invariant function $c:\Sigma\to \Bbb 
C$ such that  $\beta_s=c_s(c_s+1)$ for each $s\in \Sigma$. 
Set 
$$L=\delta_c(x)H\delta_c(x)^{-1}.$$
Then an easy computation shows 
that
$$L=\Delta-\sum_{s\in\Sigma}{2c_s\over\alpha_s(x)}\partial_{\alpha_s}$$
where, for a vector $y\in\frak h$,  as usual the symbol $\partial_y$ 
denotes the partial derivative in the $y$ direction (notice that 
using the scalar product we are viewing $\alpha_s$ as a vector in 
$\frak h$ orthogonal to the hyperplane fixed by $s$). 

From now on we will work with $L$ instead of 
$H$ and study the eigenvalue problem
\begin{equation}\label{duemezzo} L\psi=\lambda\psi\end{equation}
It is clear that
$\psi$ is a solution of this equation if and only if
$\delta_c(x)^{-1}\psi$ 
is a solution of
(\ref{due0}). 

Since for any 
$s\in \Sigma$ and $f\in \Bbb C[\frak h]$  we have that $f(sx)-f(x)$ 
is divisible by $\alpha_s(x)$, we get that the operator  
$${1\over 
\alpha_s(x)}(s-1)\in \cal A$$ maps $ \Bbb C[\frak h]$ to itself.

\begin{Definition} {\rm Given  $y\in \frak h$, we define the Dunkl 
operator
$D_y$ on $\Bbb C[\frak h]$ by

$$D_y:=\partial_y+\sum_{s\in\Sigma}c_s{(\alpha_s,y)\over 
\alpha_s(x)}(s-1)$$}
\end{Definition}

We have the following very important theorem. 

\begin {Theorem}\cite{Du}\label {Dunkl} 
Let $y, z\in\frak h$,  Then
$$[D_y,D_z]=0.$$\
\end {Theorem}
\demo  
See \cite{Du}, \cite{Op}. 
\QED

\begin {Proposition} (Cherednik) Let $\{y_1,\ldots y_n\}$ be 
an orthonormal basis of $\frak h$. Then we have

$$m(\sum_{i=1}^nD_{y_i}^2)=L.$$
\end{Proposition}
\demo  Observe that $m(\sum_{i=1}^nD_{y_i}^2)=\sum_{i=1}^nm(D_{y_i}^2) $, so we need to compute
$m(D_y^2)$ for
$y\in
\frak h$. We have $m(D_y^2)=m(D_ym(D_y))=m(D_y\partial_y)$. A simple computation shows that
$$D_y\partial_y=\partial_y^2+\sum_{s\in\Sigma}c_s{(\alpha_s,y)\over 
\alpha_s(x)}(\partial_y(s-1)-{2(\alpha_s,y)\over(\alpha_s,\alpha_s)}\partial_{\alpha_s}).$$
Thus
$$m(D_y^2)=\partial_y^2-2\sum_{s\in\Sigma}c_s{(\alpha_s,y)^2\over 
(\alpha_s,\alpha_s)\alpha_s(x)}\partial_{\alpha_s}$$
We get
$$m(\sum_{i=1}^nD_{y_i}^2)=\sum_i\partial_{y_i}^2-
2\sum_{s\in\Sigma}c_s{\sum_{i=1}^n(\alpha_s,y_{i})^2\over 
(\alpha_s,\alpha_s)\alpha_s(x)}\partial_{\alpha_s}=L$$ 
since $\sum_{i=1}^n(\alpha_s,y_i)^2= 
(\alpha_s,\alpha_s)$.
\QED

We 
are now ready to give the construction on quantum integrals of $L$. Consider 
the symmetric algebra
$S{\frak h}=\Bbb C[y_1,\ldots ,y_n]$ which 
we can identify, using the fact the the Dunkl operators commute, with 
the polynomial ring
$\Bbb C[D_{y_1},\ldots ,D_{y_n}]\subset  \cal 
A$.  The restriction of $m$ to 
$S{\frak h}^W$ is an algebra 
homomorphism   into the ring ${\cal D}(U)$ (and in fact into ${\cal 
D}(U/W)$). Since   $S{\frak h}^W$ is itself a polynomial 
ring 
$\Bbb C[q_1,\ldots ,q_n]$,  with $q_1,\ldots ,q_n$ of degree $d_1,\ldots ,d_n$, 
$d_i$ being the degrees of basic $W$-invariants,  we obtain a polynomial ring of 
commuting differential operators in  
${\cal D}(U)$. Given $q\in \Bbb 
C[q_1,\ldots ,q_n]$ we will denote by $L_q$ the corresponding 
differential operator.  We can assume that $q_1=\sum_{i=1}^ny_1^2$ so 
that $L=L_{q_1}$. Thus for every $q\in \Bbb C[q_1,\ldots ,q_n]$, 
$L_q$ is a quantum integral of the quantum Calogero-Moser system. In 
particular, the operators $L_{q_1},\ldots  L_{q_n}$ are $n$ 
algebraically independent pairwise commuting quantum integrals.

Now the eigenvalue 
problem (\ref{duemezzo})
may be replaced by 
$$L_p\psi=\lambda_p\psi$$
for $p\in \Bbb C[q_1,\ldots ,q_n]$, where 
the assignment $p\to \lambda_p$ is a algebra homomorphism $\Bbb 
C[q_1,\ldots ,q_n]\to \Bbb C.$

In other words, we may say that since 
$\Bbb C[q_1,\ldots ,q_n]=\Bbb C[\frak h^*/W]=
\Bbb C[\frak h/W]$, for 
every point $k\in \frak h/W$, we have the eigenvalue problem

\begin{equation}\label{dueuno}L_p\psi=p(k)\psi.\end{equation}

\begin {Proposition}  Near a generic point  
$x_0\in \frak h$, the system 
$L_p\psi=p(k)\psi$ has a space of solutions of dimension $|W|$.
\end {Proposition}

\demo
The proposition follows easily from the fact that 
the symbols of $L_{q_i}$ are $q_i(\partial)$, and that
$\Bbb C[y_1,...,y_n]$
is a free module over $\Bbb C[q_1,...,q_n]$ of rank $|W|$. 
\QED

\subsection{Additional integrals for integer $c$}

If $c_s\notin \Bbb Z$, the analysis of the 
solutions of the equations $L_p\psi=p(k)\psi$ 
is rather difficult (see \cite{HO}). 
However, in the case 
$c:W\to \Bbb Z$, the system can be simplified. Let us consider this 
case. First remark the since $\beta_s=c_s(c_s+1)$, by  
changing $c_s$ to $-1-c_s$ if necessary, we can assume that $c$ is non-negative. So 
we will assume that $c$ takes non-negative integral values and we 
will denote it by $m$.

System
(\ref{dueuno})
can be further 
simplified, if  we can find a differential operator $M$
(not a polynomial of $L_{q_1},...,L_{q_n}$) such 
that $[M,L_p]=0$ for all $p\in \Bbb C[q_1,\ldots ,q_n]$. Then the 
operator $M$ will act on the space of solutions of (\ref{dueuno}), 
hopefully with distinct eigenvalues. So, if $\mu$ is such an 
eigenvalue, the system
$$\left\{\matrix{  L_p\psi=p(k)\psi\cr 
M\psi=\mu\psi}\right.$$
will have a one dimensional space of 
solutions and we can find the unique up to scaling
 solution $\psi$ using Euler formula.

Such an $M$ 
exists if and only if $c=m$ has integer values. 
Namely, we will see that one can
extend the homomorphism $\Bbb C[q_1,\ldots ,q_n]\to {\cal D}(U)$ 
mapping $q\to L_q$ to the ring of $m$-quasi-invariants $Q_m$. 

We start by remarking that under some natural homogeneity 
assumptions, if such an extension exists,  it is unique. 
\begin {Proposition} 1) Assume  that $q\in \Bbb C[y_1,\ldots 
,y_n]$ is a homogeneous polynomial of degree $d$. If 
there exists a differential operator $M_q$ with coefficients in $\Bbb 
C(\frak h)$  of the form 
$$M_q=q(\partial_{y_1},\ldots 
\partial_{y_n})+l.o.t.$$ such that $[M_q,L]=0$, 
whose homogeneity degree is $-d$,
then $M_q$ 
is unique.

 2) Let $ \Bbb C[q_1,\ldots ,q_n]\subseteq 
B\subseteq \Bbb C[y_1,\ldots ,y_n]$ be a graded ring. Assume that we 
have a linear map $M:B\to {\cal D}(U)$ 
such that, if $q\in 
B$ is homogeneous of degree $d$,  then $[M_q,L]=0$, $M_q$ has homogeneity degree 
$-d$, and 
$$M_q=q(\partial_{y_1},\ldots ,
\partial_{y_n})+l.o.t.$$
Then $M$ is a ring homomorphism and 
$M_q=L_q$ for all $q\in \Bbb C[q_1,\ldots ,q_n].$ 
\end{Proposition}

\demo 1) If there exist two different operators $M_q$ and $M_q'$ with 
these properties, take $M_q-M_q'$. This operator has degree of 
homogeneity $-d$, but order smaller than $d$. Therefore, its
symbol $S(x,y)$ is not a polynomial. On the other hand, 
since the symbol of $L$ is $\sum
y_i^2$, we get that  $[L,M_q-M_q']=0$ implies $\{ \sum
y_i^2,S(x,y)\}=0$. 
Write $S$ in the form $K(x,y)/ H(x)$ with $K$ is a polynomial,
and $H(x)$ a homogeneous polynomial
of positive degree $t$ (we
assume that $K(x,y)$ and $H(x)$ have no common irreducible
factors). Then 
$$0=\{ \sum y_i^2,S(x,y)\}=2{\sum_{i=1}^ny_iK_{x_i}(x,y)H(x)-
\sum_{i=1}^ny_iH_{x_i}(x)K(x,y)\over H(x)^2}$$
Since $\sum_{i=1}^nx_iH_{x_i}(x)=tH(x)$, we have that $\sum_{i=1}^ny_iH_{x_i}(x)K(x,y)\neq 0$. So 
$H(x)$ must divide this polynomial and, by our assumptions, this implies that it must divide the 
polynomial $\sum_{i=1}^ny_iH_{x_i}(x)$ whose degree in
$x$ of
 is $t-1$. This is a contradiction.

2) Let $q,p\in B$ be two homogeneous elements. Then $M_qM_p$ and 
$M_{pq}$ both satisfy the same homogeneity assumptions. Hence they 
are equal by 1). 

Finally if $q\in 
\Bbb C[q_1,\ldots ,q_n]$, both $M_q$ and $L_q$ satisfy the same 
homogeneity assumptions. Hence they are equal by 1).
\QED

The required 
extension  to the ring of $m$-quasi-invariants is then provided by the following:

\begin {Theorem}{\rm (\cite{CV1,CV2})} Let $c=m:\Sigma\to \Bbb Z_+$. 
The following two conditions are equivalent for a homogeneous 
polynomial $q\in \Bbb C[\frak h^*]$ of degree $d$. 

1) There exists 
a differential operator 
$$L_q=q(\partial_{y_1},\ldots 
\partial_{y_n})+l.o.t.$$
of homogeneity degree $-d$,
such that 
$[L_q,L]=0$.

2) $q$ is an $m$-quasiinvariant homogeneous of degree 
$d.$
\end{Theorem}
Using this, we can extend system 
(\ref{dueuno}) to the system

\begin{equation}\label{duedue}L_p\psi=p(k)\psi,\ \ \ \ \ \  p\in 
Q_m,\ \  k\in {\rm Spec}\ Q_m=X_m\end{equation}
(Recall that, as a 
set, $X_m=\frak h$). Near a generic point $x_0\in \frak h$, system (\ref{duedue})  has a one 
dimensional space of solutions, thus there exists a unique up to scaling
solution $\psi(k,x)$, which can be expressed in elementary
functions. This solution is called the {\it Baker-Akhiezer
  function}, and has the form
$$\psi(k,x)=P(k,x)e^{(k,x)}$$
with 
$P(k,x)$ a polynomial of the form $\delta(x)\delta(k)+l.o.t.$. 
Furthermore, it can be shown that $\psi(k,x)=\psi(x,k)$ (see \cite{CV1,CV2,FV}).

These results motivate the following terminology. 
The variety $X_m$ is called {\it the spectral 
variety} of the Calogero-Moser system for the multiplicity function $m$, and $Q_m$ 
is called {\it the spectral ring} of this system.

\subsection{An example}

\begin {Example} \label{meq1} {\rm Let $W=\Bbb Z/2$, $\frak
h=\Bbb C$, $m=1$. As we have seen, $Q_m$ has a basis given by the 
monomials $\{x^{2i}\}\cup \{x^{2i+3}\}$, $i\geq 0$.  Let us set 
for such a monomial, $L_{x^r}=L_r$, and  $\partial={d\over
  dx}$. Then we have 

$$L_0=1, \ \ L_2=\partial^2-{2\over x}\partial,\ \ L_3=\partial^3-{3\over x}\partial^2+{3\over 
x^2}\partial$$
As for the others, $L_{2j}=L_2^j$, $L_{2j+3}=L_2^jL_3.$
(Note that $L_1$ is not defined). 
The system 
(\ref {duedue}) in this case is
$$\left\{\matrix{ \psi{''}-{2\over 
x}\psi'=k^2\psi,\cr \psi^{'''}-{3\over x}\psi{''}+{3\over 
x^2}\psi'=k^3\psi}\right.$$
The solution can easily be computed by 
first differentiating the first equation and subtracting the second,
thus obtaining the new system
$$\left\{\matrix{ \psi{''}-{2\over 
x}\psi'=k^2\psi\cr  \psi{''}-({1\over 
x}+k^2x)\psi'=-k^3x\psi}\right.$$
Taking the difference, we get the 
first order equation
$$\psi'={k^2x\over kx-1}\psi$$
whose solution 
(up to constants) is given by $\psi=(kx-1)e^{kx}$.}
\end{Example}

In fact, one can easily calculate $\psi_m$ for a general $m$. 

\begin{Proposition}
$$\psi_m(k,x)=(x\partial-2m+1)(x\partial-2m-1)\cdots 
(x\partial-1)e^{kx}$$
\end{Proposition}

\demo
We could use the direct method of Example \ref{meq1}, but it is
more convenient to proceed differently. 
Namely, we have 
$$(\partial^2-{2m\over 
x}\partial)(x\partial-2m+1)=(x\partial-2m+1)(\partial^2-{2(m-1)\over x}\partial)$$
as it is easy to 
verify directly. So using induction in $m$ with base $m=0$, we get
$$(\partial^2-{2m\over 
x}\partial)\psi_m(k,x)=(x\partial-2m+1)(\partial^2-{2(m-1)\over x}\partial)\psi_{m-1}(k,x)=
k^2\psi_m(k,x),$$
and $\psi_m(k,x)$ is our solution.
\QED 

\section{Lecture 3}

\subsection{Shift operator and construction of the Baker-Akhiezer
  function}

In Lecture 2, we have introduced the Baker-Akhiezer function $\psi(k,x)$ for the operator
$$L=\Delta-\sum_{s\in\Sigma}{2c_s\over\alpha_s(x)}\partial_{\alpha_s}.$$
The way to construct $\psi(k,x)$ is via Opdam shift operator.
Given a function $m:\Sigma\to \Bbb Z_+$, Opdam showed in \cite{Op1} that there exists a unique
$W$-invariant differential operator $S_m$ of the form
$\delta_m(x)\delta_m(\partial_x)+l.o.t.$, with $\delta_m(x)=\prod_{s\in\Sigma}\alpha_s^{m_s}$
such that
$$L_qS_m=S_mq(\partial)$$
for every $q\in \Bbb C[q_1,...,q_n]$. 
From this, if we set 
$$\psi(k,x)=S_m^{(x)}e^{(k,x)},$$
we get 
\begin{equation}\label{symsys}
L_q\psi=S_mq(\partial)e^{(k,x)}=q(k)\psi,
\end{equation}
$q\in \Bbb C[q_1,...,q_n]$.

We claim that equation (\ref{symsys}) must in fact hold for all $q\in
Q_m$. Indeed, near a generic point $x$, the functions 
$\psi(wk,x)$ are obviously linearly independent and 
satisfy (\ref{symsys}) for symmetric $q$. Thus, they are a basis 
in the space of solutions (we know that this space is
$|W|$-dimensional). Consider the matrix of $L_q$ in this basis
for any $q\in Q_m$. Since $\psi(k,x)$ is a polynomial times
$e^{(k,x)}$, this matrix must be diagonal with eigenvalues
$q(k)$, as desired.  

\begin {Example} 
{\rm As we have seen in the previous section, for $W=\Bbb Z/2$ and $\frak h=\Bbb C$,
$$S_m=(x\partial-2m+1)(x\partial-2m-1)\cdots 
(x\partial-1)$$.}
\end{Example}

\subsection{Berest's formula for $L_q$}

We are now going to give an explicit construction of the operators $L_q$ for any $q\in Q_m$.

Let us identify, using our $W$-invariant scalar product, $\frak h$ with $\frak h^*$, and let us
choose a orthonormal basis $x_1,\ldots ,x_n$ in $\frak h^*$. If $x\in \frak h^*$, we will write
$D_x$ for the Dunkl operator relative to the vector in $\frak h$ corresponding to $x$ under our
identification. Thus 
$$L=\sum_{i=1}^nD_{x_i}^2$$

\begin {Proposition}{\rm (Berest \cite{Be})}
\label{treuno} If $q\in Q_m$ is a homogeneous element of
degree $d$, then
$$(ad L)^{d+1}q=0.$$
\end {Proposition}
\demo It is enough to prove that 
$$((ad L)^{d+1}q)^{(x)}\psi(k,x)=0.$$
Indeed, by the definition of $\psi(k,x)$, we get that this implies that in the ring ${\cal D}(U)$,
$((ad L)^{d+1}q)S_m=0$, so that $(ad L)^{d+1}q=0$, since ${\cal D}(U)$ is a domain.

Given $q\in Q_m$, we will denote by $L_q^{(k)}$ the operator $q(D_{k_1},\ldots ,D_{k_n})$. Notice that
since
$\psi(k,x)=\psi(x,k)$, we have that  
$L_q^{(k)}\psi=q(x)\psi$. Thus we deduce, for $p,q,r\in Q_m$,
 $$L_qr(x)L_p\psi=L_qr(x)p(k)\psi=p(k)L_qr(x)\psi=p(k)L_qL_r^{(k)}\psi=$$$$ =p(k)L_r^{(k)}L_q\psi=
p(k)L_r^{(k)}q(k)\psi$$
It follows that
$$(ad L)^{d+1}q\psi=(-1)^{d+1}(ad (\sum_{i=1}^nk_i^2))^{d+1}L_q^{(k)}\psi$$
Since $L_{q}$ is a differential operator of degree $d$, we get that $ad
(\sum_{i=1}^nk_i^2)^{d+1}L_q^{(k)}=0$, as desired.
\QED

Notice now that the operator $(adL)^dq(x)$ commutes with $L$. Its symbol is given by
$(ad \Delta)^dq(x)=2^dd!q(\partial)$. So we deduce the following:
\begin {Corollary} (Berest's formula, \cite{Be}) \label{tredue} If $q\in Q_m$ is homogeneous of
  degree $d$, then 
$$L_q={1\over 2^dd!}(ad L)^dq(x).$$
\end {Corollary}
\demo This is clear from Proposition 2.8, once we remark that $(ad L)^dq(x)$  has the required
homogeneity.\QED

We want to give a representation theoretical interpretation of what we have just seen.  Consider the
three operators 
\begin{equation}F={\sum_{i=1}^nx_i^2\over 2},\ E= -{L\over 2},H=[E,F]\label
{esselle2}\end{equation} It is easy to check that
$[H,E]=2E$, $[H,F]=-2F$.
We deduce that the elements $E,F,H$ span an  ${ s\ell}(2)$
Lie subalgebra of ${\cal D}(U)$. Thus ${ s\ell}(2)$ acts by
conjugation on ${\cal D}(U)$. 
We can then reformulate Proposition~\ref {treuno} as follows:
\begin {Proposition}\label {trequattro} 
Any polynomial $q\in Q_m$ of degree $d$ is a lowest
weight vector for the
${ s\ell}(2)$-action of weight $-d$ and generates a finite dimensional module (necessarily of
dimension
$d+1$) for which $L_q$ is a highest weight vector.
\end {Proposition}
\demo An easy direct computation shows that
$$H=[E,F]=-\sum_{i=1}^nx_i{\partial\over\partial x_i}+C$$
where $C$ is a constant. Thus if $q$ is homogeneous of degree $d$,  we have
$[H,L_q]=dL_q$. 

This and the fact that $[L,L_q]=0$, implies that $L_q$ is a
highest weight vector
of weight $d$.
Also since $F$ is a polynomial, we deduce that ${\rm ad F}^{d+1}L_q=0$, so that $L_q$ generates a $d+1$
dimensional irreducible ${ s\ell}(2)$-module.\QED

One last property about these operators is given by:
\begin {Proposition}\cite{FV}\label{preserv} For any $q\in Q_m$, the operator $L_q$ preserves $Q_m$.
\end{Proposition} 

\demo
Let us begin 
by proving that $L$ preserves $Q_m$. 

Take $f\in Q_m$, so that for any $s\in\Sigma$,
$f-\thinspace ^{s}\negthinspace f=\alpha_s^{2m_s+1}t$, $t\in \Bbb C[\frak h]$.
Let us start by showing that $Lf$ is a polynomial. Clearly $Lf=\delta_*^{-1}q$, with $q\in {\Bbb C}[{\frak
h}]$, and $\delta_*=\prod_{s:m_s\ne 0}\alpha_s$. Since $L$ is $W$-invariant,
$Lf-^s\negthinspace (Lf)=L(f-\thinspace ^{s}\negthinspace f)$ is clearly divisible by
$\alpha_s^{2m_s-1}$ if $m_s>0$. In particular, it always 
is regular along the reflection hyperplane of $s$. On the other
hand, since
$Lf-\thinspace ^{s}\negthinspace (Lf)=\delta_*^{-1}(q+\thinspace ^{s}\negthinspace q)$, we deduce that
$q+\thinspace ^{s}\negthinspace q$ is divisible by
$\alpha_s$ if $m_s>0$. But then 
$q=((q+\thinspace ^{s}\negthinspace q)+(q-\thinspace ^{s}\negthinspace q))/2$ is divisible by $\alpha_s$
if $m_s>0$, hence it is divisible by $\delta_*$, so that
$Lf$ lies in $q\in {\Bbb C}[{\frak
h}]$.

 We have already remarked that  $L(f-\thinspace ^{s}\negthinspace f)$ is divisible by
$\alpha_s^{2m_s-1}$ if $m_s>0$. In fact
$$L(f-\thinspace ^{s}\negthinspace f)=(L\alpha_s^{2m+1})t+\alpha_s^{2m}\tilde t$$
$\tilde t$ being a suitable polynomial.

But since  $$L\alpha_s^{2m_s+1}=2m_s(2m_s+1)(\alpha_s,\alpha_s)\alpha_s^{2m_s-1}-2m_{s'}(2m_s+1)\sum_{s'\in
\Sigma}(\alpha_{s'},\alpha_s){\alpha_s^{2m_s}\over \alpha_{s'}}=$$$$
=-2m_{s'}(2m_s+1)\sum_{s'\in
\Sigma, s'\neq s}(\alpha_{s'},\alpha_s){{\alpha_s}^{2m}\over \alpha_{s'}},$$
we deduce that $L(f-\thinspace ^{s}\negthinspace f)$ is divisible by
$\alpha_s^{2m_s}$. On the other hand, since $L(f-\thinspace ^{s}\negthinspace f)=Lf-\thinspace
^{s}\negthinspace (Lf)$, this polynomial is either zero or it must vanish to odd order on the
reflection hyperplane of $s$. We deduce that it must be  divisible  by $\alpha_s^{2m_s+1}$, proving 
that $Lf\in Q_m$.

We now pass to a general $L_q$, $q\in Q_m$. We can assume that $q$ is homogeneous of, say, degree $d$.
By Corollary~\ref {tredue} we have that $L_q$ is a non zero multiple of $(adL)^d(q)$. Since both $q$
and
$L$ preserve $Q_m$, our claim follows.\QED

\subsection{Differential operators on $X_m$}

Now let us return to the algebra of differential operators ${\cal
  D}(X_m)$. 
Notice that ${\cal D}(X_m)$ contains two commutative subalgebras
(both 
isomorphic to $Q_m$). The first
is
$Q_m$ itself, the second is the subalgebra $Q_m^\dagger$ consisting of the differential operators of the form
$L_q$ with $q\in Q_m$. 
It is possible to show that:
\begin {Theorem}\cite{BEG} \label{tretre} ${\cal D}(X_m)$ is generated by $Q_m$ and $Q_m^\dagger$.
\end{Theorem}
 Notice that by Corollary~\ref {tredue} we in fact have that ${\cal D}(X_m)$ is generated
by $Q_m$ and by $L$.
\begin {Example}{\rm If $W=\Bbb Z/2$, $\frak h=\Bbb C$  we get that ${\cal D}(X_m)$
is generated by the operators 
$$x^2,x^{2m+1}, {d^2\over dx^2}-{2m\over x}{d\over dx}.$$}
\end{Example}
 Theorem~\ref {tretre} together with Proposition~\ref {trequattro}, imply: 
\begin {Corollary}\cite{BEG}  ${\cal D}(X_m)$ is  locally
finite dimensional under the action of the Lie algebra ${ s\ell}(2)$ defined in (\ref {esselle2}).
\end{Corollary}
This Corollary implies that our ${ s\ell}(2)$ action on ${\cal
  D}(X_m)$  can be integrated to an action of the
group $SL(2)$. In particular we have that
$$ \left(\matrix {0  & 1\cr -1 & 0}\right)q=L_q.$$
for all $q\in Q_m$. This transformation is a generalization of
the Fourier transform, 
since it reduces to the usual Fourier transform on differential
operators on $\frak h$ when $m=0$.
\begin {Example}{\rm If $W=\Bbb Z/2$, $\frak h=\Bbb C$, we get that the monomials 
$\{x^{2i}\}\cup\{x^{2i+2m+1}\}$ are (up to constants) all lowest weight vectors for the ${ s\ell}(2)$
action on ${\cal D}(X_m)$. $x^n$ has weight $-n$. We deduce that  
${\cal D}(X_m)$ is 
isomorphic as a ${ s\ell}(2)$-module to the direct sum of the irreducible representations of
dimension $n+1$ for $n$ even or $n=2(m+i)+1$, each with multiplicity one.}
\end{Example}

\subsection{The Cherednik algebra}
Let us now go back to the algebra $\cal A$  of operators on $U$ generated by ${\cal D}(U)$ 
and $W$. This algebra contains the Dunkl operators 
$$D_y:=\partial_y+\sum_{s\in\Sigma}c_s{(\alpha_s,y)\over 
\alpha_s}(s-1).$$ 

\begin {Lemma} \label {relaz} The following relations hold:
$$[x_i,x_j]=[D_{x_i},D_{x_j}]=0,\ \ \ \ \ \forall 1\leq i,j\leq n$$
$$[D_{x_i},x_j]=\delta_{i,j}+\sum_{s\in\Sigma}c_s{(x_i,\alpha_s)(x_j,\alpha_s)\over
(\alpha_s,\alpha_s)}s,\
\
\ \ \ \forall 1\leq i,j\leq n$$
$$wxw^{-1}=w(x),\ \ wD_yw^{-1}=D_{w(y)}, \ \ \ \forall w\in W, x\in \frak h^*,y\in \frak h$$
\end {Lemma}
\demo The proof is an easy computation, except the relations
$[D_{x_i},D_{x_j}]=0$, which follow from Theorem~\ref {Dunkl}.\QED

Thius lemma motivates the following definition. 

\begin {Definition}(see e.g. \cite{EG}) The Cherednik algebra $H_c$ is 
an associative algebra with generators $x_i,y_i, i=1,...,n$, and
$w\in W$, with defining relations 
$$[x_i,x_j]=[y_i,y_j]=0,\ \ \ \ \ \forall 1\leq i,j\leq n$$
$$[y_i,x_j]=\delta_{i,j}+\sum_{s\in\Sigma}c_s{(x_i,\alpha_s)(x_j,\alpha_s)\over
(\alpha_s,\alpha_s)}s,\
\
\ \ \ \forall 1\leq i,j\leq n$$
$$wxw^{-1}=w(x),\ \ wyw^{-1}=w(y), \ w\cdot w'=ww',\ \ \ \forall w,w'\in W, x\in
\frak h^*,y\in \frak h,
$$
\end{Definition}

This algebra was introduced by Cherednik as a rational limit of
his double affine Hecke algebra defined in \cite{Ch}. 
Notice that if $c=0$ then $H_c={\cal D}({\frak h})\rtimes \Bbb C[W]$.

Lemma \ref{relaz} implies that 
the algebra $H_c$ is equipped with a homomorphism $\phi:H_c\to
{\cal A}$, given by $w\to w,x_i\to x_i,y_i\to D_{x_i}$. 

Cherednik proved the following theorem.

\begin {Theorem}{\rm (Poincar\`e-Birkhoff-Witt theorem)} The multiplication map 
$$\mu:{\Bbb C}[\frak h]\otimes \Bbb C[\frak h^*]\otimes {\Bbb C}[W]$$
given by $\mu(f(x)\otimes g(y)\otimes w)=f(x)g(y)w$ is an isomorphism of vector spaces.
\end {Theorem}
\demo It is easy to see that the map $\mu$ is surjective. Thus, we only have to show that
it is injective. In other words, we need to show that 
monomials $x_1^{i_1}...x_n^{i_n}y_1^{j_1}...y_n^{j_n}w$ are
linearly independent in $H_c$. To do this, it suffices to show
that the images of these monomials under 
the homomorphism $\phi$,
i.e. $x_1^{i_1}...x_n^{i_n}D_{x_1}^{j_1}...D_{x_n}^{j_n}w$,
are linearly independent. 

Given an element $A\in \cal A$, writing $A=\sum_{w\in W}P_ww$ with $P_w\in {\cal D}(U)$ we define the
order of $A$, ord$A$, as the maximum of the orders of the $P_w$'s. Notice that
ord$AB\leq$ord$A$+ord$B$.   We now remark that for any sequence of non negative indices,
$(i_1,\ldots i_n)$, 
$$D_{x_1}^{i_1}\cdots D_{x_n}^{i_n}= \partial_{x_1}^{i_1}\cdots \partial_{x_n}^{i_n}+l.o.t.$$
        Indeed  this is true for $D_{x_i}$. We proceed by induction on $r=i_1+\cdots +i_n$. We can clearly
assume $i_1>0$, so by induction,
$$D_{x_1}^{i_1}\cdots D_{x_n}^{i_n}=(\partial_{x_1}+l.o.t.)( \partial_{x_1}^{i_1-1}\cdots
\partial_{x_n}^{i_n}+l.o.t.)=\partial_{x_1}^{i_1}\cdots \partial_{x_n}^{i_n}+l.o.t.$$

From this we deduce that for any pair of multiindices $I=(i_1,\ldots i_n)$, $J=(j_1,\ldots j_n)$, 
$w\in W$, setting $x_I={x_1}^{i_1}\cdots {x_n}^{i_n}$, $D_J=D_{x_1}^{j_1}\cdots D_{x_n}^{j_n}$,
$\partial_J=
\partial_{x_1}^{j_1}\cdots \partial_{x_n}^{j_n}$, we have
$$x_ID_Jw=x_I\partial_Jw+l.o.t.$$
Using this and the linear independence of the elements $x_I\partial_Jw$,  it is immediate to conclude
that the elements $x_ID_Jw$ are linearly independent, proving our claim.\QED

{\bf Remark 1.} We see that the homomorphism $\phi$ identifies 
$H_c$ with the subalgebra of $\cal A$ generated by $\Bbb C[\frak
h]$, the Dunkl operators $D_y$, $y\in \frak h$ and $W$. 

{\bf Remark 2.} Another way to state the PBW theorem is the
following. Let $F^\bullet$ be a filtration on $H_c$ defined by
${\rm deg}(x_i)={\rm deg}(y_i)=1$, ${\rm deg}(w)=0$. 
Then we have a natural surjective mapping
from $\Bbb C[{\frak h}\times {\frak h^*}]\rtimes W$ 
to the associated graded algebra ${\rm gr}(H_c)$.
The PBW theorem claims that this map is in fact an isomorphism.  

\subsection{The spherical subalgebra}

Let us now introduce the idempotent
$$e={1\over W}\sum_{w\in W}w=\in {\Bbb C}[W].$$
\begin {Definition} {\rm The spherical subalgebra of $H_c$ is the algebra $eH_ce$.}
\end {Definition}
Notice that $1\notin eH_ce$. On the other hand, since $ex=xe=e$
for $x\in eH_ce$,
$e$ is the unit for the spherical
subalgebra. We can embed both  $\Bbb C[\frak h^*]^W$ and $ \Bbb C[\frak h]^W$ in the spherical
subalgebra as follows. Take $f\in \Bbb C[\frak h^*]^W$ (the other case is identical) and set
$m_e(f)=fe$. Since $f$ is invariant, we have $efe=fe^2=fe=m_e(f)$, so that $m_e$ actually maps  $\Bbb
C[\frak h^*]^W$ to $eH_ce$. The injectivity is clear from the PBW-theorem. As for the fact that $m_e$
is a homomorphism, we have
$m_e(fg)=fge=fge^2=fege=m_e(f)m_e(g)$. From now on, we will consider both $\Bbb C[\frak h^*]^W$ and $
\Bbb C[\frak h]^W$ as subalgebras of the spherical subalgebra.

\subsection{Category $O$}
 We are now going to study representations  of the algebras $H_c$ and $eH_ce$. 

\begin {Definition}{\rm The category ${\cal O}(H_c)$(resp. ${\cal O}(eH_ce))$ is the full subcategory 
of the category of
$H_c$-modules (resp. $eH_ce$-modules) whose objects are the
modules
$M$  such that
 
\noindent 1) $M$ is finitely generated.

\noindent 2) For all $v\in M$, the subspace $\Bbb C[\frak h^*]^Wv\subset M$ is finite dimensional.}
\end{Definition}

We can define a functor
$$F:{\cal O}(H_c)\to {\cal O}(eH_ce)$$
by setting $F(M)=eM$. It is easy to show that $F(M)$ is a object
of ${\cal O}(eH_ce)$. 

We are now going to explain how to construct some modules in ${\cal O}(H_c)$ which, by analogy with the
case of enveloping algebras of semisimple Lie algebras, we will
call Whittaker and Verma modules.
First, take $\lambda\in{\frak h}^*$. Denote by $W_{\lambda}\subset W$ the stabilizer of $\lambda$. 
Take an irreducible $W_{\lambda}$ module $\tau$. We define a structure of $\Bbb C[\frak h^*]\rtimes \Bbb
C[W_{\lambda}]$-module on $\tau$  by 
$$(fw)v=f(\lambda)(wv), \forall v\in\tau,\ w\in W_{\lambda},\ f\in\Bbb C[\frak h^*].$$
It is easy to see that this action is well defined and we call this module $\lambda\#\tau$. We can
then consider the
$H_c$-module
$$M(\lambda,\tau)=H_c\otimes_{\Bbb C[\frak h^*]\rtimes {\Bbb C}[W_{\lambda}]}\lambda\#\tau$$
This is called a Whittaker module. In the special case
$\lambda=0$ (and hence $W_\lambda=W$), the module
$M(0,\tau)$ is called a Verma module. It is clear that
these are objects of  $\cal O$.
Notice that as $\Bbb C[\frak h]\rtimes {\Bbb C}[W]$-module,  $M(\lambda,\tau)=\Bbb C[\frak
h]\otimes_{\Bbb C}{\Bbb C}[W]\otimes_{{\Bbb C}[W_{\lambda]}}\tau$.
\begin {Example} {\rm If $\lambda=0$ and $\tau={\bf 1}$ is the trivial representation of $W$,  the
Verma module $M(0,{\bf 1})=\Bbb C[\frak h]$. The action of $\Bbb C[\frak h]$ is given by
multiplication, the one of
$\Bbb C[\frak h^*]$ is generated by the Dunkl operators and $W$ acts in the usual way.}
\end {Example}

\subsection{Generic $c$}

Opdam and Rouquier have recently studied the structure 
of the categories ${\cal O}(H_c)$, ${\cal O}(eH_ce)$, and found
that it is especially simple if $c$ is ``generic'' in a certain
sense. Namely, recall that for a $W$-invariant function 
$q: \Sigma\to \Bbb C^*$ one may define the {\it Hecke algebra}
${\rm He}_q(W)$ to be the quotient of the group algebra of the fundamental group 
of $U/W$ by the relations $(T_s-1)(T_s+q_s)=0$, where $T_s$ is
the image in $U/W$ of a small half-circle around the hyperplane
of $s$ in the counterclockwise direction. It is well known 
that ${\rm He}_q(W)$ is an algebra of dimension $|W|$, 
which coincides with $\Bbb C[W]$ if $q=1$. It is also known 
that ${\rm He}_q(W)$ is semisimple (and isomorphic to 
$\Bbb C[W]$ as an algebra) unless $q_s$ for some $s$ belongs to a finite set of
roots of unity depending on $W$. 

\begin{Definition} The function $c$ is said to be generic
if for $q=e^{2\pi ic}$, the Hecke algebra ${\rm He}_q(W)$ is
semisimple. 
\end{Definition}

In particular, any irrational $c$ is generic, and (more
importantly for us) an integer valued $c$ is generic
(since in this case $q=1$).

We can now state the following central result:
\begin{Theorem} ({\rm Opdam-Rouquier} \cite{OR}; see also
  \cite{BEG} for an exposition)
\label {semis} If $c$ is generic (in particular, if $c$ takes non negative integer values), then
the irreducible objects in
$\cal O$ are exactly the modules $M(\lambda,\tau)$. Moreover, the category $\cal O$ is semisimple.
\end{Theorem}

We also have 
\begin {Theorem}(\cite{OR})\label {equiv} 
If $c$ is generic then the functor $F$ is an equivalence of categories.
\end{Theorem}

From  Theorem~\ref {semis}  we can deduce:

\begin{Theorem} \cite{BEG} If $c$ is generic, then $H_c$ is a simple algebra.
\end{Theorem} 

In the case $c=0$, we get the simplicity of $\Bbb C[\frak
h\oplus\frak h^*]\rtimes \Bbb C[W]$, which is well known.

\subsection{The Levasseur-Stafford theorem and its
  generalization}

Let us now recall a result of Levasseur and Stafford:

\begin {Theorem}\cite{LS}
\label {invar} If $G$ is a finite group acting on a finite dimensional vector space $V$
over the complex numbers, then the ring ${\cal D}(V)^G$ is generated by the subrings ${\Bbb C}[V]^G$
and 
${\Bbb C}[V^*]^G$.
\end {Theorem}
As an example, notice that if we let $\Bbb Z/n\Bbb Z$ act on the complex line by multiplication by the n-th
roots of 1, we deduce that the operator $x{d\over dx}$ can be expressed as a non commutative polynomial
in the operators $x^n$ and ${d^n\over dx^n}$, a non-obvious
fact. We note also that this theorem has a purely ``quantum''
nature, i.e. the corresponding ``classical'' statement, saying
that the Poisson algebra $\Bbb C[V\times V^*]^G$ is generated,
as a Poisson algebra, by $\Bbb C[V]^G$ and $\Bbb C[V^*]^G$, 
is actually false, already for $V=\Bbb C$ and $G=\Bbb Z/n\Bbb Z$.

One can prove a similar result for the algebra $eH_ce$. 
Namely, recall that the algebra $eH_ce$ contains the subalgebras $\Bbb C[\frak h]^W$,
and $\Bbb C[\frak h^*]^W$.

\begin {Theorem}\cite{BEG} \label {gener}
 If $c$ is generic then the two subalgebras $\Bbb C[\frak h]^W$
and $\Bbb C[\frak h^*]^W$ generate $eH_ce$.
\end{Theorem}

Notice that if $c=0$, then $eH_0e=\cal D(\frak
h)^W$, so Theorem \ref{gener} reduces to the Levasseur-Stafford
theorem. 

{\bf Remark.} It is believed that this result holds without the
assumption of generic $c$. Moreover, it is known to be true for
all $c$ if $W$ is a Weyl group not of type $E$ and $F$, since in
this case Wallach proved that the corresponding classical
statement for Poisson algebras holds true. Nevertheless, 
the genericity assumption is needed for the proof, because,
similarly to the proof of the Levasseur-Stafford theorem, it is
based on the simplicity of $H_c$. 

\subsection{The action of the Cherednik algebra to quasi-invariants}

We now go back to the study of $Q_m$. Notice that the algebra $eH_me$ acts on $\Bbb C[\frak h]^W$, since
$e$ gives the $W$-equivariant projection of $\Bbb C[\frak h]$ onto $\Bbb C[\frak h]^W$. 
It is clear that this action is by differential operators. For instance,
the subalgebra  $\Bbb C[\frak h]^W\subset eH_me$ acts by
multiplication. Also, an element
$q\in \Bbb C[\frak h^*]^W\subset eH_me$ acts via the operator
$q(D_{x_1},\ldots D_{x_n})$. By definition this operator coincides with $L_q$ on $\Bbb C[\frak h]^W$.

The following important theorem shows that this action extends to
$Q_m$. 

\begin {Theorem}\cite{BEG} There exists a unique representation of the algebra $eH_me$ on
  $Q_m$ in which an element $q\in \Bbb C[\frak h]^W$ acts by
multiplication and an element $q\in \Bbb C[\frak h^*]^W$ by $L_q$.
\end {Theorem}

\demo Since by Proposition \ref{preserv}, $L_q$ preserve $Q_m$, 
we get a uniquely defined representation of the subalgebra of
$eH_me$ generated by $\Bbb C[\frak h]^W$ and $\Bbb C[\frak
h^*]^W$ on $Q_m$. The result now follows from
Theorem~\ref {gener}.\QED

\subsection{Proof of Theorem \ref{freeness}}

Finally we can prove Theorem \ref{freeness}.

To do this, observe that as an $eH_me$-module, 
$Q_m$ is in the category ${\cal O}(eH_me)$,
and $\Bbb C[{\frak h^*}]^W$ acts locally nilpotently 
in $Q_m$ (by degree arguments). We can now apply
Theorem~\ref{equiv} and Theorem~\ref {semis} and deduce that $Q_m$ is a direct sum of modules of the
form
$eM(0,\tau)$. As a $\Bbb
C[\frak h]\rtimes {\Bbb C}[W]$-module, $M(0,\tau)=\Bbb C[\frak h]\otimes
 \tau$. On the other hand, by Chevalley theorem,  there is an isomorphism $\Bbb C[\frak
h]\simeq \Bbb C[\frak h]^W
\otimes {\Bbb C}[W]$, commuting with the action of $W$ and $\Bbb C[\frak h]^W$. Thus we get an
isomorphisms of $\Bbb C[\frak h]^W$-modules
$$eM(0,\tau)\simeq (M(0,\tau))^W\simeq \Bbb C[\frak h]^W
\otimes ({\Bbb C}[W]\otimes  \tau )^W\simeq \Bbb C[\frak h]^W
\otimes   \tau$$
proving that $eM(0,\tau)$ and hence $Q_m$ is a free $\Bbb C[\frak h]^W$-module.\QED

\begin {Example} {\rm For $W=\Bbb Z/2$ and $\frak h=\Bbb C$, take the polynomials  $1,x^{2m+1}$.
Notice that $L(1)=L(x^{2m+1})=0$ while $s(1)=1$, $s(x^{2m+1})=-x^{2m+1}$, $s\in \Bbb Z/2$ being the
element of order two. It follows that $Q_m$ as a $eH_me$-module is the direct sum of ${\Bbb
C}[x^2]\oplus x^{2m+1}{\Bbb C}[x^2]$. These modules are
irreducible. 
Moreover, ${\Bbb C}[x^2]\simeq eM(0,{\bf
1})$,
$x^{2m+1}{\Bbb C}[x^2]\simeq eM(0,\varepsilon)$, $\varepsilon$ being the sign representation.}
\end{Example}

\subsection{Proof of Theorem \ref{Dsimple}}

Let $I$ be a nonzero two-sided ideal in ${\cal D}(X_m)$. First we claim
that $I$ nontrivially intersects $Q_m$. Indeed, otherwise let
$K\in I$ be a lowest order nonzero element in $I$. 
Since the order of $K$ is positive, 
there exists $f\in Q_m$ such that $[K,f]\ne 0$. 
Then $[K,f]\in I$ is of smaller order than $K$. Contradiction. 

Now let $f\in Q_m$ be an element of $I$. 
Then $g=\prod_{w\in W}{}^wf\in I$. 
But $g$ is $W$-invariant. This shows that the intersection 
$J$ of $I$ with the subalgebra $H_m$
in ${\cal D}(X_m)$ is nonzero. But $H_m$ is simple, so 
$J=H_m$. Hence, $1\in J\subset I$, and $I={\cal D}(X_m)$. 

\bibliographystyle{ams-alpha}

\begin{thebibliography}{A}  
\bibitem[Ar]{Ar} Arnold, V., Mathematical methods of classical
  mechanics, Graduate texts in Math., Springer Verlag, 1978. 
\bibitem[Be]{Be}
Yu. Berest, \textit{Huygens' principle and the bispectral problem}
in \textit{The Bispectral Problem}, CRM Proceedings and Lecture 
Notes, \textbf{14}, Amer. Math. Soc. 1998, pp. 11--30.
\bibitem[BEG]{BEG} Yu. Berest, P. Etingof, V. Ginzburg,
 \textit{Cherednik algebras and differential operators on
quasi-invariants}, math.QA/0111005.
\bibitem[Ch]{Ch} Cherednik, I.,  \textit{Double affine Hecke algebras, 
Knizhnik-Zamolodchikov equations, and Macdonald operators}, 
IMRN (Duke Math. J.) v.9 (1992), p.171-180. 
\bibitem[CV1]{CV1}
O.~A.~Chalykh and A.~P.~Veselov, \textit{
Commutative rings of partial differential operators and Lie algebras}, 
Comm. Math. Phys. \textbf{126}(3) (1990), 597--611. 
%
\bibitem[CV2]{CV2} 
O.~A.~Chalykh and A.~P.~Veselov, \textit{Integrability in the 
theory of Schr\"odinger operator and harmonic analysis},
Comm. Math. Phys. \textbf{152}(1) (1993), 29--40. 

\bibitem[Du]{Du}
C. F. Dunkl, \textit{Differential-difference operators and 
monodromy representations of Hecke algebras}, 
Pacific J. Math. \textbf{159}(2) (1993), 271--298.

\bibitem[Eis]{Eis} Eisenbud, D., Commutative algebra with a view 
toward algebraic geometry, Springer, New York, 1994. 

\bibitem[EG]{EG}  
P. Etingof and V. Ginzburg, \textit{Symplectic 
reflection algebras, Calogero-Moser space, and 
deformed Harish-Chandra homomorphism},
{\tt{math.AG/0011114}}, Invent. Math. (2001).
%
\bibitem[EG2]{EG2} P. Etingof and V. Ginzburg, 
\textit{On $m$-quasi-invariants of Coxeter groups},
Preprint {\tt{math.QA/0106175}}.
%
\bibitem[FV]{FV} 
M. Feigin and A. Veselov, \textit{Quasi-invariants of
Coxeter groups and $m$-harmonic polynomials}, 
Preprint 2001, {\tt{math-ph/0105014}}.
%
\bibitem[FeV]{FeV} 
G. Felder and A. Veselov, 
\textit{Action of Coxeter groups on $m$-harmonic polynomials
and KZ equations}, Preprint 2001, {\tt{QA/0108012}}.

\bibitem[HO]{HO} Heckman, G. J.; Opdam, E. M.,  \textit{Root systems and hypergeometric
functions. I}, Compositio Math. 64 (1987), no. 3, 329--352.

\bibitem[LS]{LS} Levasseur, T.; Stafford, J. T.  \textit{Invariant differential operators
and an homomorphism of Harish-Chandra}, J. Amer. Math. Soc. 8
(1995), no. 2, 365--372.

\bibitem[Op]{Op} Opdam, E. M., Lecture notes on Dunkl operators for real and
complex reflection groups, MSJ
Memoirs, 8. Mathematical Society of Japan, Tokyo, 2000.

\bibitem[Op1]{Op1} Opdam, E. M.,  \textit{Some applications of hypergeometric shift
operators}, Invent. Math. 98 (1989), no. 1, 1--18

\bibitem[OR]{OR}
E. Opdam and R. Rouquier, in preparation.

\bibitem[VdB]{VdB}
M. Van den Bergh, \textit{Differential operators on semi-invariants 
for tori and weighted projective spaces} in 
\textit{Topics in Invariant Theory}, Lecture Notes in Math. 
\textbf{1478}, Springer, Berlin, 1991, pp. 255--272. 

\end{thebibliography}

\end{document}